\documentclass[12pt,a4]{article}
\usepackage{amssymb}
\usepackage{amsfonts,amssymb}
\usepackage[]{graphicx}
\usepackage[cp1251]{inputenc}
\usepackage[english, russian]{babel}
\usepackage[a4paper, left=30mm, right=10mm, bottom=30mm, head=0mm]{geometry}
\begin{document}
 \thispagestyle{empty}

\begin{center}
 {\bf  А. Р. Миротин,  Е. Ю. Кузьменкова}
 \end{center}

\begin{center}
  ГАНКЕЛЕВЫ ОПЕРАТОРЫ НАД КОМПАКТНЫМИ АБЕЛЕВЫМИ ГРУППАМИ
\end{center}

Рассматривается два варианта обобщений операторов Ганкеля на случай линейно упорядоченных абелевых групп, даются критерии ограниченности и компактности этих операторов, в том числе в терминах функций ограниченной средней осцилляции,  доказана нефредгольмовость обобщенных операторов Ганкеля. Даны некоторые приложения к теории тёплицевых операторов на группах.

 Two variants of generalizations of Hankel operators to the case of linearly ordered abelian groups are considered, criteria of the boundedness and compactness of these operators are given, among them in terms of functions of bounded mean oscillation,  the nonfredholmness of  generalized  Hankel operators is proved. Some applications to the theory of Toeplitz operators on groups are given.

\begin{center}
{\bf 1. Введение}
\end{center}

Классические операторы Ганкеля   представляют собой один из важнейших классов операторов в пространствах голоморфных функций, имеющий интересные приложения к проблеме моментов, ортогональным полиномам, теории рациональной аппроксимации и другим важным разделам анализа, а также теории прогнозирования и теории управления (см., например, \cite{Pel}, \cite[часть В]{Nik1}, \cite[часть D, глава 5]{Nik2}).   Одно из равносильных определений ганкелева оператора состоит в том, что в некотором ортонормированном базисе он имеет (вообще говоря, бесконечную) ганкелеву матрицу, т. е. матрицу, элементы которой зависят лишь от суммы индексов. После открытия символа ганкелева оператора теория таких операторов в значительной степени свелась к изучению зависимости свойств оператора от геометрических и аналитических характеристик его символа. Таким образом, в этой теории взаимодействуют методы теории функций и функционального анализа. Не удивительно, что предпринимались многочисленные попытки     обобщения этих операторов
(см. \cite{Nak1} --  \cite{YCG} и обзор в \cite[с. 195 - 204]{Nik1}), как и  обобщения тесно связанных с ними операторов Тёплица (см. \cite{SbMath}).

В данной работе рассматривается два вида обобщений операторов Ганкеля, определенные  в пространствах  $l_2(X_+)$ и $H^2(G)$ соответственно, где $X_+$ есть положительный конус линейно упорядоченной дискретной абелевой группы $X$,   $G$  --- группа характеров группы $X$; рассмотрено также обобщение интегральных операторов Ганкеля на полуоси и их дискретного аналога. Установлены связи между ними и даны критерии их ограниченности и компактности, в том числе в терминах определенных в \cite{DM} пространств функций ограниченной средней осцилляции на группе $G$. Кроме того, доказана  нефредгольмовость обобщенных операторов Ганкеля. Основной целью являлось изучение ганкелевых операторов $H_\varphi$, действующих из $H^2(G)$ в $H^2_-(G)$  (как указано в \cite{Pel}, в классическом случае эта реализация является  наиболее важной). В случае, когда  группа $X$ содержит наименьший положительный элемент, для этих операторов  указанные вопросы решены полностью. При этом показано, что наличие наименьшего положительного элемента является и необходимым условием для существования нетривиальных компактных операторов вида $H_\varphi$.

\newpage
\begin{center}
{\bf 2. Вспомогательные сведения и результаты}
\end{center}

 Всюду ниже   $G$ есть нетривиальная связная
компактная абелева группа с  нормированной мерой Хаара $dx$ и
линейно упорядоченной группой характеров $X$,
$X_+$ --- положительный конус в  $X$. Другими
словами, в группе  $X$ выделена   подполугруппа $X_+$, содержащая
единичный характер $1$
и такая, что $X_+\cap X_+^{-1}=\{1\}$ и $X=X_+\cup X_+^{-1}$. При этом
полугруппа $X_+$ индуцирует в $X$  линейный порядок, согласованный
со структурой группы,  по правилу $\xi\leq\chi:=\chi\xi^{-1}\in
X_+$. Далее мы положим $X_-:= X_+^{-1}\setminus\{1\}(=X\setminus X_+$). Хорошо известно, что (дискретная) абелева группа $X$ может быть линейно
упорядочена тогда и только тогда, когда она не имеет кручения
(см., например, \cite{Rud}), что, в свою очередь, равносильно
тому, что её группа характеров $G$ компактна и связна \cite{Pont} (при этом линейный порядок
в $X$, вообще говоря, не единственен). В приложениях в роли  $X$ часто выступают подгруппы аддитивной группы $\mathbb{R}^n$,
наделенные дискретной топологией, так что $G$ является боровской компактификацией группы  $X$ (см, например, \cite{EMRS} и литературу там).

Рассмотрим пространство
$$l_2(X_+)=\{f:X_+\rightarrow\mathbb{C}:\sum\limits_{\chi\in X_+}| f(\chi)| ^2<\infty \}$$
(здесь каждая функция $f$ имеет не более чем счетный носитель) со
скалярным произведением
$$
\langle f,g\rangle=\sum\limits_{\chi\in X_+}f(\chi)\overline{g(\chi)}.
$$
Аналогично определяются и другие пространства $l_p(X_\pm)$. Ясно, что  система $\{1_{\{\chi\}}\}_{\chi\in
X_\pm}$ (через $1_A$ мы обозначаем индикатор множества $A$) индикаторов одноточечных подмножеств  является ортонормированным базисом пространства
$l_2(X_\pm)$.

{\bf Определение 1.} Пусть $k$ --- функция на $X_+$. Ганкелевой формой на $X_+$ c ядром $k$ называют комплексную
билинейную форму вида
$$A(a,b)=\sum\limits_{\chi,\eta\in X_+}k(\chi\eta)a(\chi)b(\eta),\eqno (1)$$
определенную первоначально на финитных  функциях  на $X_+$.

Через $\widehat{\varphi}$ (или $\mathcal{F}\varphi$)  мы будем обозначать преобразование Фурье функции $\varphi$ из  $L^1(G)$  или $L^2(G)$, т. е. (при соответствующей интерпретации интеграла)
$$
\widehat{\varphi}(\xi)=\int\limits_G\varphi(x)\overline{\xi(x)}dx,\ \xi\in X.
$$
Далее нам понадобится следующий результат, обобщающий классическую теорему Нехари \cite{Neh}.

{\bf Теорема} (Нехари-Вонг \cite{Wong}). {\it Ганкелева форма (1) на $X_+$
ограничена тогда и только тогда, когда ее ядро имеет вид
$k(\chi)=\widehat{\varphi}(\overline\chi),\ \chi\in X_+,$
 где $\varphi\in L^{\infty}(G)$.  При этом
$\|\varphi\|_{\infty}\leq M,$  где $M$ - константа
ограниченности формы $A.$
В частности,  норма формы $A$ равняется
$\|\varphi\|_{\infty}$ для некоторой функции $\varphi$, удовлетворяющей
указанным выше условиям}.

\newpage
\begin{center}
{\bf 3. Операторы Ганкеля в пространстве $l_2(X_+)$}
\end{center}

{\bf Определение 2.}  Оператор $\Gamma:l_2(X_+)\to l_2(X_+)$, определенный первоначально на финитных  функциях  на $X_+$,
называется ганкелевым (оператором Ганкеля) в $l_2(X_+)$, если
существует функция $a=a_\Gamma$ на $X_+$ такая, что для всех $\chi,\xi\in
X_+$ выполняется равенство
$$
\langle\Gamma 1_{\{\chi\}},1_{\{\xi\}}\rangle=a(\chi\xi).
$$

Классические ганкелевы операторы соответствуют случаю
 $X=\mathbb{Z}$. С точки зрения обобщений интерес
представляет уже случай  $X=\mathbb{Z}^n$ (в связи с этим отметим, что линейные порядки в $\mathbb{Z}^n$, согласованные со структурой группы, описаны в \cite{Teh}, \cite{Zajt}).
Мы дадим обобщение на случай пространства $l_2(X_+)$ классических критериев ограниченности
ганкелевых операторов. Первое из этих обобщений есть простое следствие теоремы  Нехари-Вонга и тоже может рассматриваться как абстрактная версия теоремы Нехари.

{\bf Теорема 1.} {\it Оператор Ганкеля $\Gamma$  ограничен в $l_2(X_+)$ тогда и
только тогда, когда существует функция $\psi\in L^{\infty}(G)$
такая, что $a(\chi)=\widehat{\psi}(\chi)$ для любого
$\chi\in X_+$. При этом справедливо равенство}
$$\|\Gamma\|=\inf\{\|\psi\|_{\infty}:\psi\in L^{\infty}(G),\ \widehat{\psi}(\chi)=a(\chi) \forall\chi\in X_+\}.$$

{\it Доказательство.} Первое утверждение теоремы составляет содержание леммы 6 статьи [8].
 В доказательстве этой леммы была рассмотрена билинейная форма ($f,g$ ---  финитные  функции  на $X_+$)
$$
A(f,g):=\langle\Gamma f;\bar{g}\rangle.
$$
и показано, что  форма $A$, а вместе с ней и оператор $\Gamma$, ограничена тогда и только тогда, когда существует
функция $\psi_1\in L^{\infty}(G)$ такая,что для любого $\chi\in X_+$ выполняется равенство
$a(\chi)=\widehat{\psi_1}(\bar{\chi})$. При этом $\|A\|=\|\Gamma\|$.

Для завершения доказательства  заметим, что из теоремы  Нехари-Вонга следует, что
$$
\|A\|=\inf\{\|\psi_1\|_{\infty}:a(\chi)=\widehat{\psi_1}(\bar\chi)\ \forall\chi\in X_+\},
$$
и осталось положить в этой формуле $\psi(x)=\psi_1(x^{-1})$.  Теорема
доказана.

Для формулировки еще одного критерия ограниченности ганкелевых операторов в $l_2(X_+)$  нам необходима некоторая подготовка.

{\bf Определение 3.}  Пространство Харди $H^p(G)\ (1\leq p\leq\infty)$  над
$G$ определяется следующим образом (см., например, \cite{Rud}):
$$
H^p(G)=\{f\in L^p(G):\widehat{f}(\chi)=0\ \forall\chi\in X_-\}.
$$
Обозначим через $H^2_-(G)$ ортогональное дополнение подпространства
$H^2(G)$ пространства $L^2(G)$. Тогда
$$H^2_-(G)=\{f\in L^2(G):\widehat{f}(\chi)=0\ \forall\chi\in X_+\}.$$

Ясно, что $X_+$ является ортонормированным базисом
пространства $H^2(G)$,  а  $X_-$ --
ортонормированным базисом пространства $H^2_-(G)$.
Через $P_+$ и $P_-$ мы будем обозначать ортопроекторы из  $L^2(G)$ на $H^2(G)$ и $H^2_-(G)$ соответственно.

 Перенесем теперь на функции, определенные на $G$, понятие  ограниченной средней осцилляции. Для этого напомним определение  преобразования
Гильберта на группе $G$. Мы ограничимся  случаем линейного порядка на $X$, принадлежащим С. Бохнеру и Г. Хелсону,  более общая теория
построена в \cite[глава 6]{garman} и  \cite{ijpam}.
 Для любой функции $u$ из $L^2(G,\mathbb{R})$
существует единственная функция $\widetilde{u}$ из $L^2(G,\mathbb{R})$, такая, что
$\widehat{\widetilde{u}}(1)=0$ и $u+{\it i}\widetilde{u}\in H^2(G)$.
Функция $\widetilde{u}$ называется
{\it гармонически сопряженной с $u$}.
 Линейное отображение,
получаемое в результате продолжения отображения $u\mapsto \widetilde{u}$ на
(комплексное) $L^2(G)$ по линейности, называется  {\it преобразованием Гильберта} на группе $G$. Этот оператор ограничен в $L^2(G)$.

Следующие определения мотивированы известной теоремой Ч. Феффермана \cite{F} (см. также \cite[с. 189]{Nik1}).

\textbf{Определение 4} \cite{DM}. Определим пространства $BMO(G)$ функций ограниченной средней осцилляции и $BMOA(G)$ функций ограниченной средней осцилляции аналитического типа на группе  $G$ следующим образом:
$$
BMO(G):=\{f+\widetilde{g}: f,g\in L^{\infty}(G)\},
$$
$$
BMOA(G):=BMO(G)\cap H^1(G).
$$

В \cite{DM} доказано, что эти пространства, наделенные нормой $\|\varphi\|_{BMO}=\inf\{\|f\|_\infty+\|g\|_\infty :\varphi=f+\widetilde{g}, f,g\in L^{\infty}(G)\}$, являются банаховыми.

Теперь может быть установлена следующая

{\bf Теорема 2.} {\it Оператор Ганкеля} $\Gamma$ {в $l_2(X_+)$  \it ограничен тогда и
только тогда, когда функция} $\varphi:=\sum\limits_{\chi\in
X_+}a(\chi)\chi$ {\it принадлежит} $BMOA(G)$.

{\it Доказательство.} Докажем необходимость. Пусть $\Gamma$
ограничен. Из теоремы 1 следует, что существует функция $\psi\in
L^{\infty}(G)$ такая, что $a(\chi)=\widehat{\psi}(\chi)$ для всех
$\chi\in X_+$. Тогда $\varphi=\sum_{\chi\in
X_+}\widehat{\psi}(\chi)\chi=P_+\psi$. Кроме того, поскольку
$\psi\in L^{\infty}(G)\subset L^2(G)$, то $\varphi\in L^2(G)\subset
L^1(G)$. А так как функция $\widehat{\varphi}$ сосредоточена на $X_+$, то
$\varphi\in H^1(G)$. Из  [8, лемма 2] следует, что $
i\widetilde{\psi}=2P_+\psi-\psi-\widehat{\psi}(1)$,
а потому функция $P_+\psi=\frac{1}{2}(\psi+{\it
i}\widetilde{\psi}+\widehat{\psi}(1))$ имеет вид $f+\widetilde{g}$,
где $f,g\in L^{\infty}(G)$ и, стало быть, принадлежит $BMO(G)$. Таким образом,
$\varphi\in BMOA(G)$.

Докажем достаточность. Пусть $\varphi\in BMOA(G)$, т. е.
$\varphi=f+\widetilde{g},\ f,g\in L^{\infty}(G),\varphi\in H^1(G)$. Тогда $\varphi$ принадлежит $L^2(G)$, а потому и $H^2(G)$. С учетом леммы 2 имеем
$$
\varphi=P_+f+P_+(-{\it
i}(P_+g-P_-g-\widehat{g}( 1)))=
P_+(f-{\it i}g+{\it i}\widehat{g}(1)).
$$
Функция $\psi:=f-{\it i}g+{\it i}\widehat{g}(1)$ принадлежит
$L^{\infty}(G)$ и $
\varphi=P_+\psi=\sum_{\chi\in
X_+}\widehat{\psi}(\chi)\chi,$
т. е. $a(\chi)=\widehat{\psi}(\chi)$ для всех $\chi\in X_+$.  Ограниченность $\Gamma$ теперь следует из теоремы 1.

 Ганкелевы операторы в $l_2(X_+)$ выделяются с помощью некоторого семейства коммутационных соотношений.

 \textbf{Определение 5}. Определим в пространстве $l_2(X_+)$ оператор сдвига на характер $\chi$ из $X_+$  равенствами
 \[
 {\cal S}_\chi f(\xi)=f(\chi^{-1}\xi), \mbox{ если } \chi^{-1}\xi\in X_+,\  {\cal S}_\chi f(\xi) =0 , \mbox{ если } \chi^{-1}\xi\notin X_+.
 \]
Легко проверить, что сопряженный оператор имеет вид
\[
{\cal S}_\chi^* f(\xi)=f(\chi\xi).
\]

 \textbf{Лемма 1}. \textit{Оператор $T$ в $l_2(X_+)$ будет ганкелевым тогда и только тогда,
 когда при всех  $\chi$ из $X_+$  выполняются коммутационные соотношения}
 \[
 {\cal S}_\chi^* T=T{\cal S}_\chi.
 \]

Доказательство. Если $T$  ганкелев, то
\[
\langle  {\cal S}_\chi^* T  1_{\{\xi\}}, 1_{\{\eta\}}\rangle=\langle  T  1_{\{\xi\}}, 1_{\{\chi\eta\}}\rangle=a(\xi\chi\eta).
\]
И так же легко проверяется, что $\langle T{\cal S}_\chi 1_{\{\xi\}}, 1_{\{\eta\}}\rangle=a(\chi\xi\eta)$.

Обратно, если оператор $T$ удовлетворяет указанным коммутационным соотношениям, то
\[
\langle  T  1_{\{\xi\}}, 1_{\{\eta\}}\rangle=\langle  T {\cal S}_\xi 1_{\{1\}}, 1_{\{\eta\}}\rangle=
\langle  {\cal S}_\xi^*T  1_{\{1\}}, 1_{\{\eta\}}\rangle=\langle T  1_{\{1\}}, 1_{\{\xi\eta\}}\rangle.
\]
Лемма доказана.

Рассмотрим вопрос о фредгольмовости ганкелевых операторов в $l_2(X_+)$. Напомним, что ограниченный оператор $T$ в гильбертовом пространстве $H$ называется \textit{фредгольмовым слева}, если его образ в алгебре Калкина $\mathfrak{C}(H)$  обратим слева, другими словами, если найдутся такой ограниченный оператор $L$ и компактный оператор $K$ в $H$, что $LT=I+K$. Это равносильно тому, что оператор $T$  имеет замкнутый образ и конечномерное ядро (см., например, \cite[с. 207 -- 208]{conw}).

\textbf{Теорема 3.} \textit{Ограниченный ганкелев оператор  в $l_2(X_+)$ не фредгольмов слева.}

Доказательство.  Если допустить, что некоторый ограниченный ганкелев оператор $\Gamma$ в $l_2(X_+)$  фредгольмов слева, то для некоторого ограниченного оператора $L$ и компактного $K$ будем иметь при всех  $\chi\in X_+$ в силу леммы 1
$$
L{\cal S}_\chi^* \Gamma 1_{\{1\}}=L\Gamma{\cal S}_\chi 1_{\{1\}}=L\Gamma 1_{\{\chi\}}=1_{\{\chi\}}+K1_{\{\chi\}}.\eqno(2)
$$
Мы утверждаем, что $\lim_{\chi\in X_+}\|{\cal S}_\chi^* \Gamma 1_{\{1\}}\|=0$. В самом деле, по теореме 1 найдется такая $\psi\in L^{\infty}(G)$, что
$$
\Gamma 1_{\{1\}}(\xi)=\langle \Gamma 1_{\{1\}}, 1_{\{\xi\}} \rangle=a(\xi)=\widehat\psi(\xi)\ (\xi\in X_+).
$$
Поэтому ${\cal S}_\chi^* \Gamma 1_{\{1\}}(\xi)=\widehat\psi(\chi\xi)\ (\chi, \xi\in X_+)$.
Далее, в силу теоремы Планшереля $\sum_{\eta\in X}|\widehat\psi(\eta)|^2=\|\psi\|_{L^2(G)}^2$, а потому  ряд  $\sum_{\eta\in X_+}|\widehat\psi(\eta)|^2$ сходится. Следовательно, для любого $\varepsilon>0$ найдется такое конечное множество $F_\varepsilon\subset X_+$, что   $\sum_{\eta\notin F_\varepsilon}|\widehat\psi(\eta)|^2<\varepsilon$.  Но тогда при $\chi> \max F_\varepsilon, \xi\in X_+$ будем иметь $\chi\xi> \max F_\varepsilon$ и, стало быть, при указанных $\chi$ справедливо неравенство
$$
\|{\cal S}_\chi^*\Gamma 1_{\{1\}}\|^2=\sum\limits_{\xi\in X_+}|\widehat\psi(\chi\xi)|^2<\varepsilon,
$$
что и доказывает наше утверждение.

С другой стороны, направленность $\left(1_{\{\chi\}}\right)_{\chi\in X_+}$ стремится к нулю   в слабой топологии пространства  $l_2(X_+)$, так как   с учетом сходимости ряда $\sum_{\eta\in X_+}|f(\eta)|^2$,  для любого $f\in l_2(X_+)$ имеем
$$
\lim_{\chi\in X_+}\langle 1_{\{\chi\}}, f\rangle=\lim_{\chi\in X_+}f(\chi)=0.
$$
 Следовательно, $\lim_{\chi\in X_+}\|K 1_{\{\chi\}}\|= 0$ в силу компактности $K$.  Так как $\|1_{\{\chi\}}\|=1$, то, переходя в (2)  к $\lim_{\chi\in X_+}$, приходим к противоречию.

\textbf{Следствие 1. } \textit{Существенный спектр Фредгольма оператора Ганкеля $\Gamma$ в $l_2(X_+)$ содержит нуль.}

\textbf{Следствие 2. } \textit{Ограниченный ганкелев оператор  в $l_2(X_+)$ не имеет ограниченного левого обратного.}

\textbf{Замечание 1.} Из леммы 3, доказываемой ниже, и теоремы 1 из \cite{YCG}  следует, что ненулевой компактный ганкелев оператор  в $l_2(X_+)$ существует тогда и только тогда, когда  группа $X$ содержит наименьший положительный элемент.

\begin{center}
{\bf 4. Ограниченность операторов Ганкеля, действующих из $H^2(G)$ в $ H^2_-(G)$}
\end{center}

В этом разделе мы рассмотрим еще одну версию обобщения операторов Ганкеля на случай компактных связных абелевых групп.

{\bf Определение 6} \cite{Dyba}. Пусть $\varphi\in L^{2}(G)$.
Оператором Ганкеля (ганкелевым оператором) с символом $\varphi$ назовем оператор
$H_{\varphi}:H^2(G)\rightarrow H^2_-(G)$, определяемый первоначально на  пространстве ${\rm Pol}_+(G)$ тригонометрических полиномов аналитического типа (линейных комбинациях характеров из $X_+$) равенством
$$
H_{\varphi}=P_-M_{\varphi},
$$
где  $M_{\varphi}:f\mapsto\varphi f$ --- оператор умножения на $\varphi$.

 Очевидно, что
оператор $
H_{\varphi}$ с символом $\varphi\in L^{\infty}(G)$  ограничен, причем $\|H_{\varphi}\|\leq \|\varphi\|_\infty$.  В \cite{Dyba} показано, что ограниченный оператор
$A:H^2(G)\rightarrow H^2_-(G)$ будет ганкелевым с символом $\varphi\in L^{\infty}(G)$, если и только если
функция  $\langle A\chi, \bar\xi\rangle_{L^2(G)}$ зависит лишь от $\chi\xi\
(\chi, \xi\in X_+)$  (\cite[ следствие теоремы 2.1]{Dyba}). Кроме того, $H_{\varphi+\psi}=H_{\varphi}$, если и только если $\psi\in H^{\infty}(G)$.

Известно (см., например, \cite[гл. 1, \S 8]{Pel}), что интегральные операторы Ганкеля на полуоси унитарно эквивалентны операторам Ганкеля $H_{\varphi}$ на окружности. Опишем общую схему, под которую подпадает упомянутая выше классическая ситуация. Пусть $R$ --- локально компактная абелева группа, $\widehat{R}$ --- двойственная ей группа. Будем предполагать, что $R=R_+\cup  R_-$, где множества $R_+$ и $R_-$ $\mu$-измеримы и $\mu(R_+\cap R_-)=0$ ($\mu$ --- мера Хаара группы $R$). Тогда $L^2(R)=L^2(R_+)\oplus L^2(R_-)$ (мы рассматриваем $L^2(R_\pm)$ как подпространства пространства  $L^2(R)$, считая функции из этих пространств равными нулю на дополнении к $R_\pm$).

{\bf Определение 7.} Пусть $\nu$ ---   регулярная борелевская мера на  $R$. Оператор
$$
{\cal G}_\nu f:=1_{R_-}(\nu\ast f),
$$
 определенный первоначально лишь на непрерывных функциях на $R_+$ с компактным носителем, о котором мы будем предполагать, что он действует из $L^2(R_+)$ в $L^2(R_-)$ (так будет, например, в случае конечной меры), будем называть интегральным оператором Ганкеля на группе $R$.

Ниже для функций $g:R\to \mathbb{C}$ положено $g_0(r):=g(r^{-1})$.

\textbf{Теорема 4.} \textit{Пусть $\nu$ ---  конечная  мера, и при сделанных выше предположениях  о группе $R$ существует унитарный оператор ${\cal U}:L^2(G)\to L^2(\widehat R)$ такой, что}

1) \textit{произведение  ${\cal U}\chi\cdot  \overline{{\cal U}\bar\xi}$ зависит только от} $\chi\xi\ (\chi, \xi\in X_+)$;

2) ${\cal U}( H^2(G))={\cal F}^{-1}(L^2(R_+))$,

\noindent
\textit{где ${\cal F}$  обозначает $L^2$-преобразование Фурье на $\widehat R$.}

\textit{Тогда оператор ${\cal G}_\nu$ унитарно эквивалентен некоторому ограниченному оператору Ганкеля  в $H^2(G)$.
Точнее, оператор $A:H^2(G)\to H^2_-(G)$, определяемый формулой
$$
A=({\cal F}{\cal U})^{-1}{\cal G}_\nu{\cal F}{\cal U},
$$
равняется некоторому оператору Ганкеля $H_{\varphi}$   с символом} $\varphi\in L^{\infty}(G)$.

{\it Если, кроме того, группа $X$
содержит наименьший положительный элемент (т.~е. наименьший из элементов, удовлетворяющих условию $\chi>1$), то
$$
\|{\cal G}_\nu\|=\inf\{\|\psi_1\|_{\infty}:\psi_1\in L^{\infty}(G),\widehat{\psi_1}|X_-=q\},\eqno(3)
$$
где} $q(\xi)=\int_R({\cal F}({\cal U}\overline{\xi}\cdot\overline{{\cal U}1}))_0 d\nu$,  $\xi\in X_-.$

Доказательство. В силу известных свойств свертки оператор ${\cal G}_\nu$ корректно определен и ограничен. Поэтому оператор $A$ тоже ограничен. Кроме того, он отображает $H^2(G)$ в $H^2_-(G)$.  В самом деле, унитарный оператор ${\cal F}{\cal U}$ отображает $H^2(G)$ на $L^2(R_+)$ в силу 2). Следовательно, обратный оператор отображает $L^2(R_-)=L^2(R)\ominus L^2(R_+)$ на $H^2_-(G)=L^2(G)\ominus H^2(G)$, откуда и следует наше утверждение. С целью применить упомянутое выше следствие теоремы 2.1 из \cite{Dyba},  покажем, что функция $\langle A\chi, \bar\xi\rangle_{L^2(G)}$ зависит только от $\chi\xi\ (\chi, \xi\in X_+)$. Имеем
$$
\langle A\chi, \bar\xi\rangle_{L^2(G)}=\langle ({\cal F}{\cal U})^{*}{\cal G}_\nu{\cal F}{\cal U}\chi, \bar\xi\rangle_{L^2(G)}=\langle{\cal G}_q{\cal F}{\cal U}\chi, {\cal F}{\cal U}\bar\xi\rangle_{L^2(R)}.
$$
Если мы положим для краткости  $x={\cal U}\chi, y={\cal U}\bar\xi$ , то $x, y\in {\cal F}^{-1}(L^2(R_+))$ в силу 2). По теореме Фубини при $f\in L^2(R_+),  g\in L^2(R_-)$ справедливо равенство
$$
\langle{\cal G}_\nu f,g\rangle_{L^2(R)}=\int\limits_Rf_0\ast \overline{g}\,d\nu.
$$
 Поэтому
$$
\langle A\chi, \bar\xi\rangle_{L^2(G)}=\langle{\cal G}_\nu{\cal F}x, {\cal F}y\rangle_{L^2(R)}=\int\limits_R({\cal F}x)_0\ast \overline{{\cal F}y}d\nu
$$
$$
=
\int\limits_R({\cal F} x)_0\ast ({\cal F}\bar y)_0d\nu=\int\limits_R({\cal F}(x\cdot\bar y))_0d\nu,\eqno(4)
$$
и для доказательства совпадения оператора $A$ с некоторым оператором  $H_\varphi$ осталось заметить, что $x\cdot \bar y={\cal U}\chi\cdot \overline{{\cal U}\bar\xi}$  зависит только от $\chi\xi$ в силу 1).

Пусть теперь группа $X$
содержит наименьший положительный элемент. Тогда в силу теоремы 6, доказываемой ниже,
$$
\|{\cal G}_\nu\|=\inf\{\|\psi_1\|_{\infty}:\psi_1\in L^{\infty}(G),\widehat{\psi_1}(\xi)=\widehat{\varphi}(\xi)\ \forall\xi\in X_-\}.
$$
С другой стороны, легко проверить, что $\langle H_\varphi \chi,\overline\xi\rangle=\widehat\varphi(\overline{\chi\xi})\ (\chi, \xi\in X_+)$. Полагая здесь и в (4) $\xi=1$, получаем, что $\widehat\varphi(\overline{\chi})=\int_R({\cal F}({\cal U}\chi\cdot\overline{{\cal U}1}))_0 d\nu\ (\chi\in X_+)$, откуда и следует формула (3). Теорема доказана.

\textbf{Пример 1} (интегральные операторы Ганкеля на полуоси). В классическом случае $G=\mathbb{T}$ (окружность),  $R=\mathbb{R}, R_\pm=\mathbb{R}_\pm$  Фурье-образ ${\cal F}^{-1}(L^2(\mathbb{R}_+))$ есть пространство  Харди $H^2$ в верхней полуплоскости (теорема Пэли-Винера), и в качестве унитарного отображения пространства $H^2$ в единичном круге (которое, как известно,  естественным образом отождествляется с пространством $H^2(\mathbb{T})$) на последнее пространство можно взять  ${\cal U}f(t)=\pi^{-1/2}f\circ \tau(t)/(t+i)$, где $\tau(w)=(w-i)/(w+i)$ --- конформное отображение верхней полуплоскости на единичный круг (см., например, \cite{Pel}, гл. 1, \S 8).

\textbf{Пример 2} ("интегральные"\ операторы Ганкеля на дискретной группе). Пусть $R=X, R_{\pm}=X_{\pm}$. Тогда $\widehat R=G, L^2(R_{\pm})=l_2(X_{\pm})$. Поскольку мера в теореме 4 предполагается ограниченной, можно считать, что $\nu\in l_1(X)$ (и даже  $l_1(X_-)$). Из теоремы Планшереля следует, что $\mathcal{F}^{-1}(l_2(X_+))=H^2(G)$. Очевидно, что если $\mathcal{U}$ --- единичный оператор в $L^2(G)$, то все условия теоремы 4 выполнены. Следовательно, оператор $\mathcal{G}_\nu:l_2(X_+)\rightarrow l_2(X_-)$, который сейчас имеет вид
  $$
  \mathcal{G}_\nu f(\xi)=\sum\limits_{\chi\in X_+}\nu(\xi\chi^{-1})f(\chi)\  (\xi\in X_-),
  $$
  унитарно эквивалентен некоторому оператору Ганкеля $H_{\varphi}$   с символом  $\varphi\in L^{\infty}(G)$. При этом $H_{\varphi}={\cal F}^{-1}{\cal G}_\nu{\cal F}$. Кроме того, если группа $X$
содержит наименьший положительный элемент, то из (3) следует, что
  $$
  \|{\cal G}_\nu\|=\inf\{\|\psi_1\|_{\infty}:\psi_1\in L^{\infty}(G),\widehat{\psi_1}|X_-=\nu\},
  $$
  поскольку, как легко проверить, в данном случае $q=\nu$.  (Операторы, рассмотренные в этом примере, родственны операторам Винера-Хопфа, изучавшимся в \cite{Adukov}).

  Условие  $\nu\in l_1(X)$  в примере 2 может быть ослаблено до $\nu\in l_2(X)$, и при этом мы получим другую реализацию операторов $H_\varphi$ с символами из $L^2(G)$.

  {\bf Теорема 5.} \textit{Оператор $H_\varphi$ с символами $\varphi\in L^2(G)$ унитарно эквивалентен оператору ${\cal G}_\nu$, где $\nu=\widehat \varphi$, и наоборот, для любого  $\nu\in l_2(X)$ оператор ${\cal G}_\nu$ унитарно эквивалентен некоторому оператору $H_\varphi$ с символами $\varphi=\mathcal{F}^{-1}\nu$.
  }

Доказательство. Теорема будет доказана, если мы покажем, что для любого аналитического полинома $q\in {\rm Pol}_+(G)$ и любой функции $\varphi\in L^2(G)$ справедливо равенство
$$
H_\varphi q=\mathcal{F}^{-1}{\cal G}_{\widehat \varphi}\mathcal{F}q.
$$
Для этого заметим, что для любой  функции $g\in l_2(X)$ выполняется равенство
$$
\varphi \mathcal{F}^{-1} g= \mathcal{F}^{-1}(\widehat{\varphi}\ast g).
$$
Кроме того, в силу теоремы Планшереля
$$
\mathcal{F}P_-\mathcal{F}^{-1}g=1_{X_-}g\ (g\in l_2(X)).
$$
Следовательно, при  $g\in l_2(X)$ имеем
$$
\mathcal{F}H_\varphi \mathcal{F}^{-1}g=\mathcal{F}P_-(\varphi \mathcal{F}^{-1} g)=\mathcal{F}P_-\mathcal{F}^{-1}(\widehat{\varphi}\ast g)={\cal G}_{\widehat \varphi}g.
$$
Для завершения доказательства осталось положить в последнем равенстве $q=\mathcal{F}^{-1}g$
и заметить, что $q$ пробегает ${\rm Pol}_+(G)$, когда $g$ пробегает множество всех финитных функций из $l_2(X_+)$.

\textbf{Следствие 3.} \textit{Для любой функции $\varphi\in L^2(G)$, удовлетворяющей условию $\widehat{\varphi}\in l_1(X)$, оператор $H_\varphi$ ограничен.}

В самом деле,  оператор ${\cal G}_{\widehat \varphi}$ ограничен в силу неравенства $\|\widehat \varphi\ast f\|_2\leq \|\widehat \varphi\|_1\|f\|_2$.

\textbf{Следствие 4.}  \textit{Оператор ${\cal G}_\nu$ ограничен, если $\nu\in l_2(X),\ \mathcal{F}^{-1}\nu\in L^\infty(G)$.}

Это сразу следует из того, что  оператор $H_\varphi$ ограничен, если $\varphi\in L^\infty(G)$.

Следующая теорема дает критерий ограниченности оператора $H_{\varphi}$ при условии, что группа
 $X$
содержит наименьший положительный элемент (такие группы описаны в \cite{Adukov}). Для ее доказательства нам понадобится лемма, последнее утверждение которой имеет и самостоятельный интерес.

{\bf Лемма 2.} [8, лемма 5] \textit{ Пусть $X$
содержит наименьший положительный элемент $\chi_1$.}

а) \textit{Справедливо равенство
$X_-=X_+^{-1}\chi_1^{-1}$.
Следовательно, отображения
$$
i:\{1_{\{\chi\}}\}_{\chi\in X_+}\rightarrow X_+:1_{\{\chi\}}\mapsto\chi,\
j:\{1_{\{\chi\}}\}_{\chi\in X_+}\rightarrow
X_-:1_{\{\chi\}}\mapsto{\chi}^{-1}\chi_1^{-1}$$
 продолжаются единственным образом  до изоморфизмов гильбертовых пространств}
$$
i:l_2(X_+)\rightarrow H^2(G),\
j:l_2(X_+)\rightarrow H_-^2(G).
$$

б) \textit{Оператор $H_\varphi$ с символом $\varphi\in L^\infty (G)$ унитарно эквивалентен ограниченному  ганкелеву оператору в $l_2(X_+)$ (и потому не является фредгольмовым слева) и обратно.}

\textbf{Следствие 5.} \textit{Пусть группа $X$
содержит наименьший положительный элемент. В условиях теоремы 4 оператор $\mathcal{G}_\nu$  не  фредгольмов слева. Оператор $\mathcal{G}_\nu$ на группе $X$ при $\nu\in l_2(X_-),\ \mathcal{F}^{-1}\nu\in L^\infty(G)$ не  фредгольмов слева.}

{\bf Теорема 6.} {\it Пусть} $\varphi\in L^2(G)$ {\it и $X$
содержит наименьший положительный элемент. Следующие утверждения равносильны:}

 1) {\it оператор $H_{\varphi}$ ограничен;}

 2) {\it существует функция} $\psi_1\in L^{\infty}(G)$ {\it такая,
что} $\widehat{\psi_1}|X_-=\widehat{\varphi}|X_-$;

 3) $P_-\varphi\in BMO(G).$

{\it Кроме того, если выполнено одно из условий 1) --- 3), то}
$$
\|H_{\varphi}\|=\inf\{\|\psi_1\|_{\infty}:\psi_1\in L^{\infty}(G),\widehat{\psi_1}(\xi)=\widehat{\varphi}(\xi)\ \forall\xi\in X_-\}.
$$

 Доказательство. Воспользуемся обозначениями и результатами, содержащимися в  лемме 2 и ее доказательстве (см. [8, лемма 5]). Докажем равносильность  утверждений 1) и 2).
Оператор $H_{\varphi}$
ограничен тогда и только тогда, когда  оператор $\Gamma=j^{-1}H_{\varphi}i$ ограничен.  В доказательстве утверждения б) из [8, лемма 5] показано, что последний  оператор является  ганкелевым в $l_2(X_+)$. Значит, по теореме 1 оператор $\Gamma$ ограничен тогда и только тогда,
когда существует функция $\psi\in L^{\infty}(G)$ такая, что
$\widehat{\psi}(\chi)=a_\Gamma(\chi)=\widehat{\varphi}(\overline{\chi}\chi_1^{-1})$ для
всех $\chi\in X_+$  (последнее равенство установлено в доказательстве леммы 5 из [8]). Положим $\psi_1(x)=\chi_1^{-1}(x)\psi(x^{-1})$. Тогда
$\widehat{\psi_1}(\bar{\chi}\chi_1^{-1})=\widehat{\psi}(\chi)$ для любого
$\chi\in X_+$, а потому
$\widehat{\psi_1}(\bar{\chi}\chi_1^{-1})=\widehat{\varphi}(\bar{\chi}\chi_1^{-1})$. Таким образом, равенство
$\widehat{\psi}(\chi)=\widehat{\varphi}(\overline{\chi}\chi_1^{-1})$ для
всех $\chi\in X_+$ эквивалентно равенству $\widehat\psi_1(\xi)=\widehat{\varphi}(\xi)$ для всех $\xi\in X_-$, и
равносильность утверждений 1) и 2) доказана.

Наконец, заметим, что
$\|\psi\|_{\infty}=\|\psi_1\|_{\infty}$, где, как выше,   $\psi_1(x)=\chi_1^{-1}(x)\psi(x^{-1})$.
 Тогда в силу теоремы 1
$$
\|H_{\varphi}\|=\|\Gamma\|=\inf\{\|\psi\|_{\infty}:\widehat{\psi}(\chi)=a(\chi)\
\forall\chi\in X_+\}
$$
$$
=\inf\{\|\psi_1\|_{\infty}:\widehat{\psi_1}(\xi)=\widehat{\varphi}(\xi)\
\forall\xi\in X_-\}.
$$

Докажем теперь равносильность утверждений 2) и 3). Пусть выполнено
утверждение 2). Из леммы 1 следует, что
$P_-\psi_1=(\psi_1-\widehat{\psi_1}(1))/2+
(\widetilde{i\psi_1/2})\in BMO(G).$
Это влечет утверждение 3), поскольку
$$
P_-\psi_1=\sum\limits_{\xi\in
X_-}\widehat{\psi_1}(\xi)\xi=\sum\limits_{\xi\in
X_-}\widehat{\varphi}(\xi)\xi=P_-\varphi.
$$

Пусть теперь выполнено  3). Тогда $P_-\varphi=f+\widetilde{g}$, где $f, g\in L^\infty(G)$. Снова применяя лемму 1, имеем
$$
P_-\varphi=P_-f+P_-(-
i(P_+g-P_-g-\widehat{g}(1)))=
P_-(f+ ig+
i\widehat{g}(1)).
$$
Рассмотрим функцию $\psi_1=f+ ig+
i\widehat{g}(1)\in L^{\infty}(G)$. Так как
$P_-\psi_1=P_-\varphi$,  то $\widehat{\psi_1}(\xi)=\widehat{\varphi}(\xi)$ для всех $\xi\in X_-$, что и завершает  доказательство.

\textbf{Следствие 6.} {\it Если  $X$
содержит наименьший положительный элемент, то ганкелев оператор $H_\varphi$ с символом $\varphi\in L^2(G)$ ограничен тогда и только тогда, когда он совпадает с ганкелевым оператором $H_{\psi_1}$ с символом  $\psi_1\in L^\infty(G)$.}

Доказательство.  Это следует из  теоремы 6 и  равенства  $H_\varphi =H_{P_-\varphi}$.

\begin{center}
{\bf 5. Компактность операторов Ганкеля, действующих из $H^2(G)$ в $ H^2_-(G)$}
\end{center}

Далее для функции $f:G\to \mathbb{C}$ положим
$$
(Jf)(x)=f(x^{-1}).
$$

Следующий класс операторов рассматривался в  \cite{YCG} (там они обозначались $H_\varphi$).

\textbf{Определение 8.} Для   $\psi\in L^\infty (G)$ определим оператор $\Gamma_\psi:H^2(G)\to H^2(G)$ формулой
$$
\Gamma_\psi(f)=P_+J(\psi f).
$$

\textbf{Лемма 3}. \textit{Оператор $\Gamma_\psi\ (\psi\in L^\infty (G))$ унитарно эквивалентен некоторому ограниченному  ганкелеву оператору в $l_2(X_+)$ и наоборот, ограниченный  ганкелев оператор $\Gamma$ в $l_2(X_+)$ унитарно эквивалентен некоторому оператору $\Gamma_{J\psi}$, где $\psi\in L^\infty (G)$,   причем $a_\Gamma=\widehat{\psi}|X_+$.}

Доказательство.
Пусть $U:H^2(G)\to l_2(X_+)$ --- такой изоморфизм гильбертовых пространств, что $U\chi=1_{\{\chi\}}$  при всех $\chi\in X_+$. Если мы положим $\Gamma:=U\Gamma_\psi\ U^{-1}$, то  получим ограниченный ганкелев оператор в  $l_2(X_+)$, поскольку
\[
\langle \Gamma 1_{\{\chi\}}, 1_{\{\xi\}}\rangle=\langle U\Gamma_\psi\chi, U \xi\rangle=\langle P_+J(\psi\chi),\xi\rangle=\langle J(\psi\chi),\xi\rangle
$$
$$
=
\int\limits_{G}\psi(x^{-1})\chi(x^{-1})\overline{\xi(x)}dx=\widehat{\psi}(\overline{\chi\xi}).
\]

 Обратно, пусть $\Gamma$ --- ограниченный ганкелев оператор в $l_2(X_+),\ T:=U^{-1}\Gamma U$.
 Тогда
 \[
 \langle T \chi,  \xi\rangle=\langle U\Gamma 1_{\{\chi\}}, UU^{-1}1_{\{\xi\}}\rangle=\langle \Gamma 1_{\{\chi\}}, 1_{\{\xi\}}\rangle=a_\Gamma(\chi\xi).
 \]
 По теореме 1 существует такая функция $\psi\in L^\infty (G)$, что $a_\Gamma=\widehat{\psi}|X_+$. Поэтому, с учетом предыдущих выкладок,
 $$
 \langle T \chi,  \xi\rangle=\widehat{\psi}(\chi\xi)=\langle \Gamma_{J\psi}\chi, \xi\rangle
 $$
 для любых $\chi, \xi\in X_+$, а потому $T=\Gamma_{J\psi}$.

В  \cite{YCG}  показано (теорема 1.2), что компактный ненулевой оператор вида $\Gamma_\psi$ при $\psi\in L^\infty$ существует тогда и только тогда, когда в группе $X$ имеется наименьший положительный элемент $\chi_1$. При этом оператор $\Gamma_\psi$ компактен тогда и только тогда, когда $\psi\in H^\infty(G)+C_e$, где $C_e$ --- замкнутая подалгебра алгебры $L^\infty(G)$, порожденная полиномами от $\overline{\chi_1}$. Установим аналогичные результаты для операторов  $H_{\varphi}:H^2(G)\rightarrow H^2_-(G)$. Ниже $K_1(G)$ --- замкнутая подалгебра алгебры $L^\infty(G)$, порожденная  элементом  $\overline{\chi_1}$. Ясно, что $K_1(G)=C_e\overline{\chi_1}$.

{\bf Теорема 7.} \textit{I. Компактный ненулевой оператор $H_{\varphi}:H^2(G)\rightarrow H^2_-(G)$  с символом $\varphi\in L^\infty(G)$ существует тогда и только тогда, когда в группе $X$ имеется наименьший положительный элемент}.

\textit{II. Пусть $X$
содержит наименьший положительный элемент $\chi_1$, а $\varphi\in L^2(G)$. Следующие утверждения равносильны}:

1) {\it оператор $H_{\varphi}$ компактен;}

2) $P_-\varphi\in K_1(G).$

Доказательство. I. Необходимость. Используя метод из \cite{YCG}, можно показать, что если компактный ненулевой оператор $H_{\varphi}$  с символом $\varphi\in L^\infty(G)$ существует, то существует
ганкелев оператор с аналогичными свойствами и символом $\chi\in X_-$.  Приведем это рассуждение для полноты изложения. Пусть $\tau(x)f(y):=f(xy)\ (x\in G)$ --- оператор сдвига в $L^2(G)$ (и $H^2(G)$). Так как $\tau(x)\chi=\chi(x)\chi \ (\chi\in X)$, то $\tau(x)$ коммутирует с $P_\pm$; в частности, $\tau(x)H_\varphi\tau(x^{-1})=H_{\tau(x^{-1})\varphi}$, а потому последний оператор компактен. Далее, так как $\varphi\notin H^\infty(G)$, то $\widehat\varphi(\chi)\neq 0$ для некоторого характера $\chi\in X_-$. В силу непрерывной (и линейной) зависимости оператора  $H_{\varphi}$ от $\varphi\in L^\infty(G)$, компактный оператор
 $$
 \int\limits_G\chi(x)H_{\tau(x^{-1})\varphi}dx
 $$
(интеграл Бохнера существует) равен $H_{\chi\ast\varphi}=\widehat\varphi(\chi)H_\chi$, откуда и вытекает доказываемое утверждение. Далее, легко проверить, что компактный оператор $S_{\overline\chi} H_\chi$ в $H^2(G)$ является проектором на подпространство, порожденное множеством характеров $[1,\chi^{-1})$ (для этого достаточно рассмотреть его значения на базисе $X_+$ пространства $H^2(G)$), а потому это множество конечно, что и доказывает необходимость.

Достаточность следует из того, что если  в группе $X$ имеется наименьший положительный элемент $\chi_1$, то образом оператора $H_{\overline{\chi_1}}$  является   одномерное подпространство $\mathbb{C}\cdot \overline{\chi_1}\subset H^2_-(G)$.

II.  Если выполнено 1), то $H_{\varphi}$ ограничен, и в силу утверждения 2) теоремы 6  найдется функция $\psi_1\in L^{\infty}(G)$  такая,
что $\widehat{\psi_1}|X_-=\widehat{\varphi}|X_-$, т.~е. $P_-\varphi=P_-{\psi_1}$. Следовательно, $H_{\psi_1}=H_{\varphi}$. В начале доказательства теоремы 6 показано, что оператор $H_{\psi_1}$  унитарно эквивалентен  некоторому ограниченному ганкелеву оператору $\Gamma$ в $l_2(X_+)$, и при этом $a_\Gamma(\chi)=\widehat{\psi_1}(\overline{\chi}\chi_1^{-1})$ для
всех $\chi\in X_+$. В свою очередь, в силу леммы 3 оператор $\Gamma$ унитарно эквивалентен некоторому оператору $\Gamma_{J\psi}$, где $\psi\in L^\infty (G)$, причем $a_\Gamma=\widehat{\psi}|X_+$. Значит, оператор $H_{\psi_1}$  унитарно эквивалентен  оператору $\Gamma_{J\psi}$, причем
$$
\widehat{\psi}(\chi)=\widehat{\psi_1}(\overline{\chi}\chi_1^{-1}) \mbox{ для
всех } \chi\in X_+.\eqno(5)
$$
В силу компактности оператора $\Gamma_{J\psi}$, по теореме 1.2 из  \cite{YCG} имеем $J\psi\in H^\infty(G)+C_e$,  а потому  $P_-(J\psi)\in C_e$.  Далее, из (5) вытекает, что $\widehat{\psi_1}(\xi\chi_1^{-1})= \widehat{J\psi}(\xi)$ при $\xi\in X_-\cup\{1\}$. Поскольку $X_-\chi_1^{-1}=X_-\setminus\{\chi_1^{-1}\}$ (см. лемму 2), то отсюда следует, что
\[
P_-(J\psi)=\sum_{\xi\in X_-}\widehat{\psi_1}(\xi\chi_1^{-1})\xi=\sum_{\xi\in X_-\chi_1^{-1}}\widehat{\psi_1}(\zeta)\zeta\chi_1
\]
\[
=\left(\sum_{\xi\in X_-\setminus\{\chi_1^{-1}\}}\widehat{\psi_1}(\zeta)\zeta\right)\chi_1=(P_-\psi_1-\widehat{\psi_1}(\chi_1^{-1})\chi_1^{-1})\chi_1.
\]
Следовательно, $P_-\psi_1=(P_-(J\psi)+\widehat{\psi_1}(\chi_1^{-1}))\chi_1^{-1}\in K_1(G)$.

   Если выполнено 2), т.~е. $\varphi_-:=P_-\varphi\in K_1(G)$, то оператор $\Gamma_{\varphi_-}$ компактен. При доказательстве импликации $1 \Rightarrow 2$ было показано, что оператор $H_{\varphi_-}$ унитарно эквивалентен оператору $\Gamma_{J\psi}$, где $\psi\in L^\infty(G), \widehat{\psi}(\chi)=\widehat{\varphi_-}(\overline{\chi}\chi_1^{-1})$ для
всех $\chi\in X_+$  (см. (5)). Следовательно,
$$
P_-(J\psi)=\sum_{\chi\in X_+\setminus\{1\}}\widehat{\psi}(\chi)\bar{\chi}=\sum_{\chi\in X_+\setminus\{1\}}\widehat{\varphi_-}(\bar{\chi}\chi_1^{-1})\bar{\chi}
$$
$$
=\sum_{\xi\in X_-\chi_1^{-1}}\widehat{\varphi_-}(\xi)\xi\chi_1=\sum_{\xi\in X_-\setminus\{\chi_1^{-1}\}}\widehat{\varphi_-}(\xi)\xi\chi_1=\varphi_-\chi_1-\widehat{\varphi_-}(\chi_1).
$$
Поэтому $P_-(J\psi)\in C_e$. Кроме того, $P_+(J\psi)\in H^2(G)\cap L^\infty(G)=H^\infty(G)$, так как $P_+(J\psi)=J\psi-P_-(J\psi) \in L^\infty(G)$. Снова применяя теорему 1.2 из   \cite{YCG}, получаем, что оператор $\Gamma_{J\psi}$, а вместе с ним и $H_\varphi$, компактен, что и требовалось доказать.

\textbf{Замечание 2.}  Для случая $\varphi\in L^2(G)$ вопрос о справедливости  утверждения,
аналогичного  утверждению I теоремы 7, остается открытым.

\textbf{Следствие 7.} I. \textit{Пусть $\nu\in l_2(X), \mathcal{F}^{-1}\nu\in L^\infty(G)$.  Отличный от нуля компактный оператор $\mathcal{G}_\nu$ на группе $X$ существует тогда и только тогда, когда $X$
содержит наименьший положительный элемент.}

II. \textit{Пусть группа $X$
содержит наименьший положительный элемент, $\nu\in l_2(X)$. Оператор $\mathcal{G}_\nu$ компактен тогда и только тогда, когда}
$$
\sum\limits_{\xi\in X_-}\nu(\xi)\xi\in K_1(G).
$$

В самом деле, оператор $\mathcal{G}_\nu$ унитарно эквивалентен некоторому оператору  $H_\varphi$,  причем по теореме 5  $\widehat{\varphi}=\nu$, а потому $P_-\varphi=\sum_{\xi\in X_-}\nu(\xi)\xi$.

{\bf Следствие 8}.  ,
 \textit{Пусть $\nu\in l_2(\Bbb{Z}^d),\ d\geq 1$. Тогда оператор  $\mathcal{G}_\nu$ на группе $\Bbb{Z}^d_{{\rm lex}}$ (индекс ${\rm lex}$ указывает, что порядок лексикографический) компактен, если и только если  функция $\nu$ сосредоточена на множестве  $\{(0,\dots,0,n_{d}): n_{d}<0\}$, и}
$$
\sum_{n_d<0}\nu(0,\dots,0,n_d)t_d^{n_d}\in P_-(C(\mathbb{T})).
$$

Доказательство. Пусть $G=\Bbb{T}^d$. Тогда $X=\Bbb{Z}^d_{{\rm lex}}$. В  этом случае в $X$ существует наименьший положительный элемент $(0,\dots,0,1)$, положительный конус  имеет вид
 $$
   \Bbb{Z}^d_{{\rm lex}+}=\{n\in \mathbb{Z}^d| n_1>0\}\sqcup\{n\in \mathbb{Z}^d| n_1=0, n_2>0\}
   \sqcup
$$
$$
   \cdots\sqcup \{n\in \mathbb{Z}^d| n_1= n_2=\ldots=n_{d-1}=0,
   n_d>0\}\sqcup\{0\},
 $$
а отрицательный --- вид
 $$
   \Bbb{Z}^d_{{\rm lex}-}=\{n\in \mathbb{Z}^d| n_1<0\}\sqcup\{n\in \mathbb{Z}^d| n_1=0, n_2<0\}
   \sqcup
$$
$$
   \cdots\sqcup \{n\in \mathbb{Z}^d| n_1= n_2=\ldots=n_{d-1}=0,
   n_d<0\}.
 $$
Если отождествить точку $n\in \Bbb{Z}^d$  с характером
$\xi_n(t)=t_1^{n_1}\dots t_d^{n_d}$ группы  $\Bbb{T}^d$, то наименьший положительный характер суть $\chi_1(t)=t_{d}$. Следовательно, замкнутая подалгебра $C_e\subset C(\Bbb{T}^d)$, порожденная полиномами от $\overline{\chi_1}$, состоит из функций на $\Bbb{T}^d$, которые не зависят от $t_1,\dots,t_{d-1}$ и (как функции от  $t_d$) принадлежит $P_-(C(\mathbb{T}))+\mathbb{C}$. Поэтому алгебра $K_1(\Bbb{T}^d)=\overline{t_d}C_e$ состоит из функций $f\in C_e$,  удовлетворяющих условию $f(0)=0$, т. е. равна $P_-(C(\mathbb{T}))$. Теперь ясно, что критерий компактности
$$
\sum\limits_{n\in  \Bbb{Z}^d_{{\rm lex}-}}\nu(n)t_1^{n_1}\dots t_d^{n_d}\in K_1(\Bbb{T}^d)
$$
из следствия 7 выполняется тогда и только тогда, когда  $\nu$ равна нулю вне множества  $\{(0,\dots,0,n_{d}): n_{d}<0\}$ и $\sum_{n_d<0}\nu(0,\dots,0,n_d)t_d^{n_d}\in P_-(C(\mathbb{T}))$.

\begin{center}
{\bf 6. Некоторые приложения}
\end{center}

Дадим приложения теоремы 7 к теории операторов Тёплица и вычислению некоторых алгебр Бургейна над группами.

Пусть $\varphi\in L^\infty(G)$.  \textit{Оператор Тёплица} $T_\varphi$ в $H^2(G)$ определяется равенством
$$
T_\varphi =P_+M_\varphi.
$$

\textbf{Следствие 9.} {\it Если функция $\varphi$ принадлежит алгебре $K_1(G)+H^\infty(G)$ и обратима в ней, то тёплицев оператор $T_\varphi$  в $H^2(G)$ фредгольмов, причем $T_{\varphi^{-1}}$ является его регуляризатором.}

Доказательство. Сопряженный к $H_\varphi$ оператор, как легко проверить, имеет вид $H_\varphi^*=P_+M_{\overline \varphi}|H^2_-(G)$. Поэтому, как и в классическом случае (см., например, \cite[лемма 4.4.3]{Nik1}), при $\varphi, \psi\in L^\infty(G)$ для полукоммутатора тёплицевых операторов справедливо равенство
$$
[T_\varphi, T_\psi):=T_\varphi T_\psi-T_{\varphi\psi}=-H^*_{\overline \varphi}H_\psi. \eqno(6)
$$
Отсюда в силу теоремы 7 будет следовать, что операторы $T_\varphi T_{\varphi^{-1}}-I$ и $T_{\varphi^{-1}}T_\varphi-I$ компактны, если показать, что для $\varphi\in L^\infty(G)$ условие 2) этой теоремы равносильно включению $\varphi\in K_1(G)+H^\infty(G)$. Докажем это. Если  $\varphi\in K_1(G)+H^\infty(G)$, то $P_-\varphi\in K_1(G)$, так как $K_1(G)\subset H^2_-(G)$. Обратно, если $P_-\varphi\in K_1(G)$, то $P_-\varphi\in L^\infty(G)$, а потому и $P_+\varphi\in L^\infty(G)$. Следовательно,  $P_+\varphi\in H^\infty(G)$, и $\varphi=P_+\varphi+P_-\varphi\in K_1(G)+H^\infty(G)$, что и требовалось.

Равенство (6) и теорема 7 влекут также следующие два утверждения:

\textbf{Следствие 10.} \textit{Пусть $\varphi, \psi\in  L^\infty(G)$. Операторы $H^*_{\overline \varphi}H_\psi$ и $[T_\varphi, T_\psi)$  компактны, если $\overline \varphi\in K_1(G)+H^\infty(G)$ или $\psi\in K_1(G)+H^\infty(G)$.}

\textbf{Следствие 11.} \textit{Пусть $\varphi\in  L^\infty(G)$. Оператор $H^*_{\overline \varphi}H_\psi$ компактен для любой $\psi\in L^\infty(G)$ тогда и только тогда, когда $\overline \varphi\in K_1(G)+H^\infty(G)$.}

\textbf{Следствие 12.} \textit{Пусть $\varphi\in  L^\infty(G)$. Следующие утверждения равносильны}:

1) \textit{коммутатор $[P_+,M_\varphi]$ компактен};

2) \textit{ коммутатор $[P_+-P_-,M_\varphi]$ компактен};

3) $\overline \varphi, \varphi\in K_1(G)+H^\infty(G)$.

Доказательство  следствия 12 не отличается от классического случая, см. \cite[следствие 2.3.3]{Nik1}.

Теорема 7 может быть использована для вычисления алгебры Бургейна над пространством $L^2(G)$, рассматриваемым как $L^\infty(G)$-модуль.
Следуя  \cite[замечание 2]{CJY}, примем следующее

\textbf{Определение 9.} Пусть $A$ --- банахова алгебра, $M$ --- банахов $A$-модуль, и  $Y$ --- замкнутое подпространство  $M$. Элемент $a\in A$  называется бургейновым, если ${\rm dist}_M(ay_n,Y)\to 0$, лишь только последовательность  $y_n$ стремится к $0$ слабо в $Y$. Множество всех бургейновых элементов обозначается $Y_b^A$ и называется алгеброй Бургейна.

\textbf{Следствие 13.}\textit{ Имеет место равенство}
$$
H^2(G)_b^{L^\infty(G)}=K_1(G)+H^\infty(G).
$$

Доказательство. Положим $A=L^\infty(G), M=L^2(G), Y=H^2(G)$. Поскольку
$$
 {\rm dist}_{L^2(G)}(\varphi y_n, H^2(G))=\|P_-(\varphi y_n)\|=\|H_\varphi y_n\|,
$$
то бургейновость элемента $\varphi\in L^\infty(G)$ равносильна полной непрерывности, т. е. компактности, оператора $H_\varphi$. Следствие 13 теперь вытекает из теоремы 7 (как уже упоминалось, для $\varphi\in L^\infty(G)$ условие 2) этой теоремы равносильно включению $\varphi\in K_1(G)+H^\infty(G)$).

\end{document}